\documentclass[11pt,a4paper]{amsart}
\usepackage{amsmath}

\usepackage{amssymb,latexsym}
\usepackage[dvips]{graphicx}

\newtheorem{theorem}{Theorem} [section]

\newtheorem{proposition}{Proposition} [section]

\title[First eigenvalue of the Laplacian on the Klein bottle]
{Greatest least eigenvalue of the Laplacian on the Klein
bottle}

\author[A. El Soufi, H. Giacomini and M. Jazar]{Ahmad El Soufi,
Hector Giacomini and Mustapha Jazar}
\address{Ahmad El Soufi and Hector Giacomini: Universit\'e
Fran\c{c}ois Rabelais de Tours, Laboratoire de Math\'ematiques et
Physique Th\'eorique, UMR-CNRS 6083, Parc de Grandmont, 37200
Tours, France.}
\address{Mustapha Jazar: Lebanese University, Mathematics
Department, P.O. Box 155-012, Beirut, Lebanon.}

\email{ahmad.elsoufi@univ-tours.fr,
hector.giacomini@univ-tours.fr, mjazar@ul.edu.lb}



\begin{document}

\begin{abstract}
We prove the following conjecture recently formulated by Jakobson,
Nadirashvili and Polterovich \cite{JNP}: For any Riemannian metric
$g$ on the Klein bottle $\mathbb{K}$ one has
$$\lambda_1 (\mathbb{K}, g) A (\mathbb{K}, g)\le
12 \pi E(2\sqrt 2/3),$$ where $\lambda_1(\mathbb{K},g)$ and
$A(\mathbb{K},g)$ stand for the least positive eigenvalue of the
Laplacian and the area of $(\mathbb{K},g)$, respectively, and $E$
is the complete elliptic integral of the second kind. Moreover,
the equality is uniquely achieved, up to dilatations, by the
metric
$$g_0= {9+ (1+8\cos ^2v)^2\over 1+8\cos^2v} \left(du^2 + {dv^2\over 1+8\cos ^2v}\right),$$
with $0\le u,v <\pi$. The proof of this theorem leads us to study
 a Hamiltonian dynamical system which turns out to be completely
integrable by quadratures.
\end{abstract}

\maketitle

\section{Introduction and statement of main results}

Among all the possible Riemannian metrics on a compact
differentiable manifold $M$, the most interesting ones are those
which extremize a given Riemannian invariant. In particular, many
recent works have been devoted to the metrics which maximize the
fundamental eigenvalue $\lambda_1(M,g)$ of the Laplace-Beltrami
operator $\Delta_g$ under various constraints (see, for instance,
\cite{ BLY, H, LY, N, P}). Notice that, since $\lambda_1(M,g)$ is
not invariant under scaling ($\lambda_1(M, kg) =
k^{-1}\lambda_1(M,g)$), such constraints are necessary.

In \cite{YY}, Yang and Yau proved that, on any compact orientable
surface $M$, the first eigenvalue $\lambda_1(M,g)$ is uniformly
bounded over the set of Riemannian metrics of fixed area. More
precisely, one has, for any Riemannian metric $g$ on $M$,
$$\lambda_1(M,g) A(M,g) \le 8\pi (\hbox {genus} (M) +1),$$
where $A(M,g)$ stands for the Riemannian area of $(M,g)$ (see
\cite{EI1} for an improvement of the upper bound). In the
non-orientable case, the following upper bound follows from Li and
Yau's work \cite {LY}: $\lambda_1(M,g) A(M,g) \le 24\pi (\hbox
{genus} (M) +1)$. On the other hand, if the dimension of $M$ is
greater than 2, then $\lambda_1(M,g)$ is never bounded above over
the set of Riemannian metrics of fixed volume, see \cite{CD}.

Hence, one obtains a relevant topological invariant of surfaces by
setting, for any compact 2-dimensional manifold $M$,
$$ \Lambda (M) = \sup_g \lambda_1(M,g) A(M,g)=\sup_{g\in {\mathcal R}(M)}\lambda_1(M,g),$$
where ${\mathcal R}(M)$ denotes the set of Riemannian metrics of
area 1 on $M$. The natural questions related to this invariant
are~:

\begin{enumerate}
\item How does $ \Lambda (M)$ depend on (the genus of) $M$?

\item Can one determine $ \Lambda (M)$ for some $M$?

\item Is the supremum $\Lambda (M)$ achieved, and, if so, by what
metrics?
\end{enumerate}

Concerning the first question, it follows from \cite{CE} that $
\Lambda (M)$ is an increasing function of the genus with a linear
growth rate. Explicit answers to questions (2) and (3) are only
known for the three following surfaces~: the sphere ${\mathbb
S}^2$, the real projective plane ${\mathbb R}P^2$ and the Torus
${\mathbb T}^2$. Indeed, according to the results of Hersch
\cite{H}, Li and Yau \cite{LY} and Nadirashvili \cite{N}, one has
$$ \Lambda ({\mathbb S}^2) = \lambda_1 ({\mathbb S}^2, g_{{\mathbb S}^2}) A
({\mathbb S}^2, g_{{\mathbb S}^2}) = 8\pi,$$ where $g_{{\mathbb
S}^2}$ is the standard metric of ${\mathbb S}^2$ (see \cite{H}),
$$ \Lambda ({\mathbb R}P^2)
= \lambda_1 ({\mathbb R}P^2, g_{{\mathbb R}P^2}) A({\mathbb R}P^2,
g_{{\mathbb R}P^2} ) = 12\pi,$$ where $g_{{\mathbb R}P^2} $ is the
standard metric of ${\mathbb R}P^2$ (see \cite{LY}), and
$$ \Lambda ({\mathbb T}^2) = \lambda_1 ({\mathbb T}^2, g_{eq}) A ({\mathbb T}^2, g_{eq}) =
{8\pi^2\over {\sqrt{3}}},$$ where $g_{eq}$ is the flat metric on
${\mathbb T}^2 \simeq {\mathbb R}^2 / \Gamma_{eq}$ corresponding
to the equilateral lattice $ \Gamma_{eq}={\mathbb Z}(1,0) \oplus
{\mathbb Z}({1\over 2},{\sqrt{3}\over 2}) $ (see \cite{N}).
Moreover, on each one of these three surfaces, the maximizing
metric is unique, up to a dilatation.

\medskip

{\it What about the Klein bottle $\mathbb{K}$?}

\medskip
Nadirashvili \cite{N} proved that the supremum $\Lambda
(\mathbb{K})$ is necessarily achieved by a regular (real analytic)
Riemannian metric. Recently, Jakobson, Nadirashvili and
Polterovich \cite{JNP} conjectured that the exact value of
$\Lambda (\mathbb{K})$ is given by
$$\Lambda (\mathbb{K}) = 12 \pi E(2\sqrt 2/3)\simeq13.365\,\pi,$$
where $E(2\sqrt 2/3)$ is the complete elliptic integral of the
second kind evaluated at $\frac{2\sqrt2}3$. They also conjectured
that this value is uniquely achieved, up to a dilatation, by the
metric of revolution
$$g_0= {9+ (1+8\cos ^2v)^2\over 1+8\cos ^2v} \left(du^2 + {dv^2\over 1+8\cos ^2v}\right),$$
with $0\le u,v <\pi$.

\medskip

The main purpose of this paper is to prove this conjecture.
Indeed, we will prove the following

\begin{theorem}\label{main}

For any Riemannian metric $g$ on $\mathbb{K}$ one has
$$\lambda_1 (\mathbb{K}, g) A (\mathbb{K}, g)\le
\lambda_1 (\mathbb{K}, g_{0}) A (\mathbb{K}, g_{0})= 12 \pi
E(2\sqrt 2/3).$$ Moreover, the equality holds for a metric $g$ if,
and only if, $g$ is homothetic to $g_0$.
\end{theorem}

As noticed in \cite {JNP}, if we denote by $({\mathbb T}^2, \bar
g_0)$ the double cover of $ (\mathbb{K}, g_{0})$, then $\bar g_0$
is nothing but the Riemannian metric induced on ${\mathbb T}^2$
from the bipolar surface of Lawson's minimal torus $\tau_{3,1}$
defined as the image in ${\mathbb S}^3$ of the map
$$(u, v)\mapsto ( \cos v \exp (3iu),  \sin v \exp (iu)).$$

It is worth noticing that the metric $g_0$ does not maximize the
systole functional $g\mapsto \hbox{sys}(g)$ (where $\hbox{sys}(g)$
denotes the length of the shortest noncontractible loop) over the
set of metrics of fixed area on the Klein bottle (see \cite {B}),
while on ${\mathbb R}P^2$ and ${\mathbb T}^2$, the functionals
$\lambda_1$ and sys are maximized by the same Riemannian metrics.

The proof of Theorem \ref{main} relies on the characterization of
critical metrics of the functional $\lambda_1$ in terms of minimal
immersions into spheres by the first eigenfunctions. Note that, in
spite of the non-differentiability of this functional with respect
to metric deformations, a natural notion of criticality can be
introduced (see \cite{EI3}). Moreover, the criticality of a metric
$g$ for $\lambda_1$ with respect to area preserving deformations
is characterized by the existence of a family $h_1, \cdots, h_d$
of first eigenfunctions of $\Delta_g$ satisfying $\sum_{i\le d}
dh_i\otimes dh_i =g$. This last condition actually means that the
map $(h_1, \cdots, h_d)$ is an isometric immersion from $(M,g)$
into ${\mathbb R}^d$ whose image is a minimal immersed submanifold
of a sphere. Theorem \ref{main} then follows from Nadirashvili's
existence result and the following

\begin{theorem}\label{th1}

The Riemannian metric $g_0$ is, up to a dilatation, the only
critical metric of the functional $\lambda_1$ under area
preserving deformations of metrics on $\mathbb{K}$.
\end{theorem}
Equivalently, the metric $g_0$ is, up to a dilatation, the only
metric on $\mathbb{K}$ such that $(\mathbb{K}, g_{0}) $ admits a
minimal isometric immersion into a sphere by its first
eigenfunctions.

In \cite{EI2}, Ilias and the first author gave a necessary
condition of symmetry for a Riemannian metric to admit isometric
immersions into spheres by the first eigenfunctions. On the Klein
bottle, this condition amounts to the invariance of the metric
under the natural ${\mathbb S}^1$-action on $\mathbb{K}$. Taking
into account this symmetry property and the fact that any metric
$g$ is conformally equivalent to a flat one, for which the
eigenvalues and the eigenfunctions of the Laplacian are explicitly
known, it is of course expected that the existence problem of
minimal isometric immersions into spheres by the first
eigenfunctions reduces to a second order system of ODEs (see
Proposition \ref{prop}). Actually, the substantial part of this
paper is devoted to the study of  the following second order
nonlinear system:
\begin{equation}\label{eq:1}
\left\{\begin{array}{lcl}
\displaystyle{\varphi_1''= (1-2\varphi_1^2 -8\varphi_2^2) \varphi_1}, \\
\\
\displaystyle{\varphi_2'' = (4-2\varphi_1^2 -8\varphi_2^2)
\varphi_2},
\end{array}\right.
\end{equation}
for which we look for periodic solutions satisfying
\begin{equation}\label{eq:3}
\left\{\begin{array}{l}
\varphi_1 \hbox{ is odd and has exactly two zeros in a period,} \\
\varphi_2 \hbox{ is even and positive everywhere;}
\end{array}
\right.
\end{equation}
and the initial conditions
\begin{equation}\label{eq:2}
\left\{\begin{array}{l}
\displaystyle{\varphi_1(0)=\varphi_2'(0)=0\; } \hbox{(from parity conditions),}\\
\displaystyle{\varphi_2(0)={1\over 2} \varphi_1'(0)=:p \in (0,1 ]
}.
\end{array}
\right.
\end{equation}

In \cite{JNP}, Jakobson, Nadirashvili and Polterovich proved that
the initial value $p=\varphi_2(0)=\sqrt{3/8}$ corresponds to a
periodic solution of (\ref{eq:1})-(\ref{eq:2}) satisfying
(\ref{eq:3}). Based on numerical evidence, they conjectured that
this value of $p$ is the only one corresponding to a periodic
solution satisfying (\ref{eq:3}). As mentioned by them, a
computer-assisted proof of this conjecture is extremely difficult,
due to the lack of stability of the system.

In Section 3, we provide a complete analytic study of System
(\ref{eq:1}). First, we show that this system admits two
independent first integrals (one of them has been already found in
\cite{JNP}). Using a suitable linear change of variables, we show
that the system becomes Hamiltonian and, hence, integrable. The
general theory of integrable Hamiltonian systems tells us that
bounded orbits correspond to periodic or quasi-periodic solutions
(see \cite{Arnold}). However, to distinguish periodic solutions
 from non-periodic ones is not easy in general. Fortunately, our
 first integrals turn out to be
quadratic in the momenta, which enables us to apply the classical
Bertrand-Darboux-Whittaker Theorem and, therefore, to completely
decouple the system by means of a parabolic type change of
coordinates $(\varphi_1,\varphi_2)\mapsto(u,v)$.
 We show that, for any
$p\ne\sqrt{3}/2$, the solutions $u$ and $v$ of the decoupled
system are periodic. The couple $(u,v)$ is then periodic if and
only if the periods of $u$ and $v$ are commensurable. We express
the periods of $u$ and $v$ in terms of hyper-elliptic integrals
and study their ratio as a function of $p$. The following fact
(Proposition \ref{pro1}) gives an idea about the complexity of the
situation: there exists a countable dense subset
$\mathcal{P}\subset(0,\sqrt{3}/2)$ such that the solution of
(\ref{eq:1})-(\ref{eq:2}) corresponding to $p\in(0,\sqrt3/2)$ is
periodic if and only if $p\in\mathcal{P}$.

In conclusion, we show that the solution associated with
$p=\sqrt{3/8}$ is the only periodic one to satisfy Condition
(\ref{eq:3}).

\section{Preliminaries: reduction of the problem}

Let us first recall that, if $g_\varepsilon$ is a smooth
deformation of a metric $g$ on a compact surface $M$, then the
function $\varepsilon \mapsto \lambda_1(M,g_\varepsilon)$, which
is continuous but not necessarily differentiable at $\varepsilon
=0$, always admits left and right derivatives at $\varepsilon =0$
with
$${d\over d\varepsilon } \lambda_1(M,g_\varepsilon) \big|_{\varepsilon=0^+}\leq
{d\over d\varepsilon } \lambda_1(M,g_\varepsilon)
\big|_{\varepsilon=0^-}$$ (see \cite{EI3, EI4} for details). The
metric $g$ is then said to be critical for the functional
$\lambda_1$ under area preserving deformations if, for any
deformation $g_\varepsilon$ with $g_0=g$ and $A(M,g_\varepsilon) =
A(M,g)$, one has
$${d\over d\varepsilon } \lambda_1(M,g_\varepsilon) \big|_{\varepsilon=0^+}\leq 0\le
{d\over d\varepsilon } \lambda_1(M,g_\varepsilon)
\big|_{\varepsilon=0^-}.$$
 It is easy to check that this last condition is equivalent to:
$$\lambda_1(M,g_\varepsilon)\leq\lambda_1(M,g)+o(\varepsilon)
\mbox{ as } \varepsilon\rightarrow 0.$$
 Of course, the metric $g$ is critical for the functional $\lambda_1$
under area preserving deformations if, and only if, it is critical
for the functional $ \lambda_1\cdot A$.

Following \cite{EI3} and \cite{N}, a necessary and sufficient
condition for the metric $g$ to be critical for $\lambda_1$ under
area-preserving metric deformations is that there exists a family
$h_1, \cdots, h_d$ of first eigenfunctions of $\Delta_g$
satisfying
\begin{equation}\label{A4}
\sum_{i\le d} dh_i\otimes dh_i = g,
\end{equation}
which means that the map $h=(h_1, \cdots, h_d)$ is an isometric
immersion from $(M,g)$ to ${\mathbb R}^d$. Since $h_1, \cdots,
h_d$ are eigenfunctions of $\Delta_g$, the image of $h$ is a
minimal immersed submanifold of the Euclidean sphere ${\mathbb
S}^{d-1}\left(\sqrt{2\over \lambda_1 (M,g)}\right)$ of radius
$\sqrt{{2}/{\lambda_1 (M,g)} }$ (Takahashi's theorem \cite{T}). In
particular, we have
\begin{equation}\label{A5}
\sum_{i\le d} h_i^2={2\over \lambda_1 (M,g)}.
\end{equation}

In \cite{EI2}, Ilias and the first author have studied conformal
properties of Riemannian manifolds $(M,g)$ admitting such minimal
isometric immersions into spheres. It follows from their results
that, if $g$ is a critical metric of $\lambda_1$ under area
preserving deformations, then
\begin{itemize}
\item[(i)] $g$ is, up to a dilatation, the unique critical metric
in its conformal class, \item[(ii)] $g$ maximizes the restriction
of $\lambda_1$ to the set of metrics conformal to $g$ and having
the same volume, \item[(iii)]the isometry group of $(M,g)$
contains the isometry groups of all the metrics $g'$ conformal to
$g$.
\end{itemize}

For any positive real number $a$, we denote by $\Gamma_{a}$ the
rectangular lattice of ${\mathbb R}^2$ generated by the vectors
$(2\pi, 0)$ and $(0,a)$ and by $ \tilde{g}_a$ the flat Riemannian
metric of the torus $ {\mathbb T}^2_a \simeq {\mathbb R}^2/
\Gamma_{a}$ associated with the rectangular lattice $ \Gamma_{a}$.
The Klein bottle $\mathbb{K}$ is then diffeomorphic to the
quotient of $ {\mathbb T}^2_a $ by the involution $s~: (x,y)
\mapsto (x+\pi, -y)$. We denote by $g_a$ the flat metric induced
on $\mathbb{K}$ by such a diffeomorphism. It is well known that
any Riemannian metric on $\mathbb{K}$ is conformally equivalent to
one of the flat metrics $g_a$.

Let $g=fg_a$ be a Riemannian metric on $\mathbb{K}$. From the
property (iii) above, if $g$ is a critical metric of $\lambda_1$
under area preserving deformations, then $\hbox {Isom}
(\mathbb{K},g_a)\subset \hbox {Isom} (\mathbb{K},g)$, which
implies that the function $f$ is invariant under the ${\mathbb
S}^1$-action $(x,y)\mapsto (x+t,y)$, $t\in [0,\pi]$, on
$\mathbb{K}$, and then, $f$ (or its lift to ${\mathbb R}^2$) does
not depend on the variable $x$.

\begin{proposition}\label{prop}
Let $a$ be a positive real number and $f$ a positive periodic
function of period $a$. The following assertions are equivalent
\begin{itemize}

\item[(I)] The Riemannian metric $g=f(y)g_a$ on $\mathbb{K}$ is a
critical metric of $\lambda_1$ under area preserving deformations.

\item[(II)] The function $f$ is proportional to $\varphi_1^2 +4
\varphi_2^2$, where $\varphi_1$ and $\varphi_2$ are two periodic
functions of period $a$ satisfying the following conditions:

\begin{enumerate}
\item[(a)] $(\varphi_1, \varphi_2)$ is a solution of the equations
\begin{equation*}
\left\{\begin{array}{lcl}
\displaystyle{\varphi_1'' = (1-2\varphi_1^2 -8\varphi_2^2) \varphi_1}, \\
\\
\displaystyle{\varphi_2''= (4-2\varphi_1^2
-8\varphi_2^2)\varphi_2};
\end{array}\right.
\end{equation*}
\item[(b)] $\varphi_1$ is odd, $\varphi_2$ is even and
$\varphi'_1(0)=2\varphi_2(0)$; \item[(c)] $\varphi_1$ admits two
zeros in a period and $\varphi_2$ is positive everywhere;
\item[(d)] $\varphi_1^2 +\varphi_2^2 \le 1$ and the equality holds
at exactly two points in a period.

\end{enumerate}
\end{itemize}
\end{proposition}

Most of the arguments of the proof of ``(I) implies (II)'' can be
found in \cite {JNP}. For the sake of completeness, we will recall
the main steps. The proof of the converse relies on the fact that
the system (\ref{eq:1}) admits two independent first integrals.

\begin{proof} The Laplacian $\Delta_g$ associated with
the Riemannian metric $g=f(y) g_a$ on $\mathbb{K}$ can be
identified with the operator $-{1\over f(y)} \left({\partial_x ^2
}+ {\partial_y ^2 }\right)$ acting on\\ $\Gamma_a$-periodic and
$s$-invariant functions on ${\mathbb R}^2$. Using separation of
variables and Fourier expansions, one can easily show that any
eigenfunction of $\Delta_g$ is a linear combination of functions
of the form $ \varphi_k (y)\cos kx$ and $\varphi_k(y)\sin kx$,
where, $\forall k$, $\varphi_k$ is a periodic function with period
$a$ satisfying $\varphi_k(-y) = (-1)^k \varphi_k(y)$ and
$\varphi_k'' = (k^2-\lambda f)\varphi_k$. Since a first
eigenfunction always admits exactly two nodal domains, the first
eigenspace of $\Delta_g$ is spanned by
$$\left\{\varphi_0(y),\; \varphi_1(y)\cos x,\; \varphi_1(y)\sin x,\;
\varphi_2 (y)\cos 2x,\; \varphi_2(y)\sin 2x \right\},$$ where,
unless they are identically zero, $\varphi_2 $ does not vanish
while $\varphi_0$ and $\varphi_1$ admit exactly two zeros in $[0,
a)$. In particular, the multiplicity of $\lambda_1(\mathbb{K}, g)$
is at most 5.

Let us suppose that $g$ is a critical metric of $\lambda_1$ under
area preserving deformations and let $h_1, \cdots, h_d$ be a
family of first eigenfunctions satisfying the equations (\ref{A4})
and (\ref{A5}) above. Without loss of generality, we may assume
that $\lambda_1(\mathbb{K}, g) =2$ and that $h_1, \cdots, h_d$ are
linearly independent, which implies that $d\le 5$. Since $h=(h_1,
\cdots, h_d):\mathbb{K}\to {\mathbb S}^{d-1}$ is an immersion, one
has $d\ge 4$. If $d=4$, then using elementary algebraic arguments
like in the proof of Proposition 5 of \cite{MR}, one can see that
there exists an isometry $\rho \in O(4)$ such that $\rho \circ h=(
\varphi_1(y) e^{i x}, \varphi_2(y) e^{2i x})$ with $\varphi_1^2
+\varphi_2^2 =1 $ (eq. (\ref{A5})) and ${\varphi'}_1^2
+{\varphi'}_2^2 =\varphi_1^2 +4\varphi_2^2=f $ (eq. (\ref{A4}))
which is impossible since $\varphi_1^2 +\varphi_2^2 =1 $ implies
that $\varphi_1$ and $\varphi_2$ admit a common critical point.
Therefore, $d=$ multiplicity of $\lambda_1(\mathbb{K}, g) =5$ and
there exists $\rho \in O(5)$ such that $\rho \circ h=(
\varphi_0(y), \varphi_1(y) e^{i x}, \varphi_2 (y) e^{2i x})$, with
$\varphi_0^2+ \varphi_1^2 +\varphi_2^2 =1 $ and
${\varphi'}_0^2+{\varphi'}_1^2 +{\varphi'}_2^2 =\varphi_1^2
+4\varphi_2^2=f $. Since the linear components of $\rho \circ h$
are first eigenfunctions of $(\mathbb{K}, g) $, one should has,
$\forall k=0, 1, 2$, $\varphi_k'' = (k^2-\lambda_1(\mathbb{K}, g)
f)\varphi_k = (k^2-2\varphi_1^2 -8\varphi_2^2) \varphi_k$. Now, it
is immediate to check that one of the couples of functions $
\left(\pm \varphi_1, \pm \varphi_2\right)$ satisfies the
Conditions (a), $\dots$, (d) of the statement. Indeed, the parity
condition $\varphi_k(-y) = (-1)^k \varphi_k(y)$ implies that
$\varphi_1(0)=\varphi'_0(0)=\varphi'_2(0)=0$ and, then,
${\varphi'}^2_1(0)=4\varphi^2_2(0)$. Conditions (c) and (d) follow
from the fact that a first eigenfunction has exactly two nodal
domains in $\mathbb{K}$.

Conversely, let $\varphi_1$ and $\varphi_2$ be two periodic
functions of period $a$ satisfying Conditions (a), $\dots$, (d)
and let us show that the metric $g=f(y) g_a$, where $f=
\varphi_1^2 +4\varphi_2^2$, is a critical metric of $\lambda_1$
under area preserving deformations. For this, we set
$\varphi_0=\sqrt{1- \varphi_1^2 -\varphi_2^2 }$ and define the map
$h:\mathbb{K}\to {\mathbb S}^4$ by $h=( \varphi_0(y), \varphi_1(y)
e^{i x}, \varphi_2 (y) e^{2i x})$. It suffices to check that the
components of $h$ are first eigenfunctions of $\Delta_g$
satisfying (\ref{A4}).

Indeed, in the next section we will see that the second order
differential system satisfied by $\varphi_1$ and $\varphi_2$
(Condition (a)) admits the two following first integrals:
\begin{equation}\label{first0}
\left\{\begin{array}{l}
(\varphi_1^2+4\varphi_2^2)^2-\varphi_1^2-16\varphi_2^2+{\varphi_1'}^2+4{\varphi'_2}^2= C,\\ \\
12\varphi_2^2(\varphi_2^2-1)+3\varphi_1^2\varphi_2^2+
\varphi_2^2{\varphi_1'}^2-2\varphi_1{\varphi_1'}\varphi_2\varphi_2'
+(3+\varphi_1^2){\varphi_2'}^2=C,
\end{array} \right.
\end{equation}
with $C= 4\varphi_2(0)^2(4\varphi_2(0)^2 - 3)$ (note that
Condition (b) implies that $\varphi_1(0)=\varphi'_2(0)=0$).
Differentiating $\varphi_0^2+ \varphi_1^2 +\varphi_2^2 =1 $ and
using the second equation in (\ref{first0}), we get
\begin{eqnarray*}
\varphi_0^2{\varphi_0'}^2&=&\varphi_1^2{\varphi_1'}^2+
\varphi_2^2{\varphi_2'}^2+2\varphi_1\varphi_1'\varphi_2\varphi_2'\\
&=&\varphi_1^2{\varphi_1'}^2+\varphi_2^2{\varphi_2'}^2+
12\varphi_2^2(\varphi_2^2-1)+3\varphi_1^2\varphi_2^2+
\varphi_2^2{\varphi_1'}^2+(3+\varphi_1^2){\varphi_2'}^2-C\\
&=&(\varphi_1^2+\varphi_2^2){\varphi_1'}^2+
(3+\varphi_1^2+\varphi_2^2){\varphi_2'}^2 +
12\varphi_2^2(\varphi_2^2-1)+3\varphi_1^2\varphi_2^2-C\\
&=&(1-\varphi_0^2){\varphi_1'}^2+(4-\varphi_0^2){\varphi_2'}^2
+12\varphi_2^2(\varphi_2^2-1)+3\varphi_1^2\varphi_2^2-C.
\end {eqnarray*}
Therefore
\begin{eqnarray*}
\varphi_0^2\left({\varphi_0'}^2 + {\varphi_1'}^2 +
{\varphi_2'}^2 \right)&=&{\varphi_1'}^2+
4{\varphi_2'}^2+12\varphi_2^2(\varphi_2^2-1)+3\varphi_1^2\varphi_2^2-C\\
&=& \left(1-\varphi_1^2 -\varphi_2^2\right) \left(\varphi_1^2
+4\varphi_2^2\right),
\end {eqnarray*}
where the last equality follows from the first equation of
(\ref{first0}). Hence,
$$|{\partial_y h }|^2= {\varphi_0'}^2 + {\varphi_1'}^2
+{\varphi_2'}^2= \varphi_1^2 +4\varphi_2^2 =
|{\partial_x h}|^2$$
and, since $\partial_x h$ and $\partial_y h$
are orthogonal, the map $h$ is isometric, which means that
Equation (4) is satisfied.

From Condition (a) one has $\varphi_1'' = (1- 2f)\varphi_1$ and
$\varphi_2'' = (4-2f)\varphi_2$, which implies that the functions
$h_1=\varphi_1(y) \cos x$, $h_2=\varphi_1(y) \sin x$,
$h_3=\varphi_2(y) \cos 2x$ and $h_4=\varphi_2(y) \sin 2 x$ are
eigenfunctions of $\Delta_g$ associated with the eigenvalue
$\lambda =2$. Moreover, differentiating twice the identity
$\varphi_0^2+ \varphi_1^2 +\varphi_2^2 =1 $ and using Condition
(a) and the identity ${\varphi_0'}^2 + {\varphi_1'}^2
+{\varphi_2'}^2= \varphi_1^2 +4\varphi_2^2=f$, one obtains after
an elementary computation, $ \varphi_0''=-2f \varphi_0$. Hence,
all the components of $h$ are eigenfunctions of $\Delta_g$
associated with the eigenvalue $\lambda =2$. It remains to prove
that $2$ is the first positive eigenvalue of $\Delta_g$ or,
equivalently, for each $k=0, 1, 2$, the function $\varphi_k$
corresponds to the lowest positive eigenvalue of the
Sturm-Liouville problem $\varphi''=(k^2-\lambda f) \varphi$
subject to the parity condition $\varphi(-y) = (-1)^k \varphi(y)$.
As explained in the proof of Proposition 3.4.1 of \cite{JNP}, this
follows from conditions (c) and (d) giving the number of zeros of
$\varphi_k$, and the special properties of the zero sets of
solutions of Sturm-Liouville equations (oscillation theorems of
Haupt and Sturm).
\end{proof}


\section{Study of the dynamical system: proof of results}
According to Proposition \ref{prop}, one needs to deal with the
following system of second order differential equations (Condition
(a) of Prop. \ref{prop})
\begin{equation}\label{s1}\left\{\begin{array}{l}
\varphi_1''=(1-2\varphi_1^2-8\varphi_2^2)\varphi_1,\\
\varphi_2''=(4-2\varphi_1^2-8\varphi_2^2)\varphi_2,
\end{array}\right.\end{equation}
subject to the initial conditions (Condition (b) of Prop.
\ref{prop})
\begin{equation}\label{datas1}\left\{\begin{array}{l}
\varphi_1(0)=0,\; \;\; \varphi_2(0)=p,\\
\varphi_1'(0)=2p,\,  \varphi_2'(0)=0,
\end{array}\right.
\end{equation}
where $p\in (0,1]$ (Condition (d) of Prop. \ref{prop}).

Notice that the system (\ref{s1})-(\ref{datas1}) is invariant
under the transform $$(\varphi_1 (y), \varphi_2(y))\mapsto
(-\varphi_1 (-y), \varphi_2(-y)).$$ Consequently, the solution
$(\varphi_1, \varphi_2)$ of (\ref{s1})-(\ref{datas1})  is such
that $\varphi_1$ is odd and $\varphi_2$ is even.

We are looking for periodic solutions satisfying the following
condition (Condition (c) of Prop. \ref{prop}):
\begin{equation}\label{zeros}
\left\{\begin{array}{l}

\varphi_1 \mbox{ has exactly two zeros in a period,}\\
\varphi_2 \mbox{ is positive everywhere.}
\end{array}
\right.
\end{equation}

Our aim is to prove the following

\begin{theorem}\label{th2}
There exists only one periodic solution of
(\ref{s1})-(\ref{datas1}) satisfying Condition (\ref{zeros}). It
corresponds to the initial value $\varphi_2(0)=p=\sqrt{3/8}$.
\end{theorem}

In fact, this theorem follows from the qualitative behavior of
solutions, in terms of $p$, given in the following

\begin{proposition}\label{pro1} Let $(\varphi_1, \varphi_2)$ be
the solution of (\ref{s1})-(\ref{datas1}).

\begin{enumerate}

\item For all $p\in(0,1]$, $p\ne\sqrt{3}/2$, $(\varphi_1,
\varphi_2)$ is periodic or quasi-periodic.

\item For $p=\frac{\sqrt3}2$, $(\varphi_1, \varphi_2)$ tends to
the origin as $y\to\infty$ (hence, it is not periodic).

\item For all $p\in(\sqrt3/2,1]$, $\varphi_2$ vanishes at least
once in each period (of $\varphi_2$). Hence, Condition
(\ref{zeros}) is not satisfied.

\item There exists a countable dense subset
$\mathcal{P}\subset(0,\sqrt{3}/2)$, with $\sqrt{3/8}\in
\mathcal{P}$, such that the solution $(\varphi_1, \varphi_2)$
corresponding to $p\in(0,\sqrt3/2)$ is periodic if and only if
$p\in\mathcal{P}$.

\item For $p=\sqrt{3/8}$, $(\varphi_1, \varphi_2)$ satisfies
(\ref{zeros}) and, for any $p\in \mathcal{P}$, $p\neq\sqrt{3/8}$,
$\varphi_1$ admits at least 6 zeros in a period.

\end{enumerate}
\end{proposition}

The first fundamental step in the study of the system above is the
existence of the following two independent first integrals.

\subsection{First integrals}
The functions
\begin{equation}\label{H1H2}
\left\{\begin{array}{l}
H_1(\varphi_1,\varphi_2,\varphi_1',\varphi_2'):=
(\varphi_1^2+4\varphi_2^2)^2-\varphi_1^2-16\varphi_2^2+
(\varphi_1')^2+4(\varphi_2')^2,\\ \\
H_2(\varphi_1,\varphi_2,\varphi_1',\varphi_2'):=
12\varphi_2^2(\varphi_2^2-1)+3\varphi_1^2\varphi_2^2+
\varphi_2^2(\varphi_1')^2\\ \\ \hskip
5cm-2\varphi_1\varphi_1'\varphi_2\varphi_2'
+(3+\varphi_1^2)(\varphi_2')^2,
\end{array} \right.
\end{equation}
are two independent first integrals of (\ref{s1}), i.e. they
satisfy the equation
$$\varphi_1'\frac{\partial H_i}{\partial \varphi_1}+
\varphi_2'\frac{\partial H_i}{\partial \varphi_2}+
\varphi_1''\frac{\partial H_i}{\partial
\varphi_1'}+\varphi_2''\frac{\partial H_i}{\partial
\varphi_2'}\equiv0.$$ The first one, $H_1$, has been obtained by
Jakobson et al. \cite{JNP}. The orbit of a solution of (\ref{s1})
is then contained in an algebraic variety defined by
\begin{equation}
\left\{\begin{array}{l}
H_1(\varphi_1,\varphi_2,\varphi_1',\varphi_2')=K_1,\\ \\
H_2(\varphi_1,\varphi_2,\varphi_1',\varphi_2')=K_2,
\end{array} \right.
\end{equation}
where $K_1$ and $K_2$ are two constants. Taking into account the
initial conditions (\ref{datas1}), one has
$K_1=K_2=-4p^2(3-4p^2)$. In other words, the solution of
(\ref{s1})-(\ref{datas1}) is also solution of
\begin{equation}\label{first}
\left\{\begin{array}{l}
(\varphi_1^2+4\varphi_2^2)^2-\varphi_1^2-16\varphi_2^2+
(\varphi_1')^2+4(\varphi_2')^2+4p^2(3-4p^2)=0,\\ \\
12\varphi_2^2(\varphi_2^2-1)+3\varphi_1^2\varphi_2^2+
\varphi_2^2(\varphi_1')^2-2\varphi_1\varphi_1'\varphi_2\varphi_2'\\
\\ \hskip 4cm+(3+\varphi_1^2)(\varphi_2')^2+4p^2(3-4p^2)=0,
\end{array} \right.
\end{equation}
with the initial conditions
\begin{equation}\label{datafirst}
\left\{\begin{array}{l}
\varphi_1(0)=0,\\ \\
\varphi_2(0)=p.
\end{array} \right.
\end{equation}
Notice that the parameter $p$ appears in both the equations
(\ref{first}) and the initial conditions (\ref{datafirst}). The
system (\ref{first}) gives rise to a ``multi-valued" 2-dimensional
dynamical system in the following way.

\subsection{2-dimensional dynamical systems}\label{2ddynamic}
From (\ref{first}) one can extract explicit expressions of
$\varphi_1'$ and $\varphi_2'$ in terms of $\varphi_1$ and
$\varphi_2$. For instance, eliminating $\varphi_1'$, one obtains
the following fourth degree equation in $\varphi_2'$
\begin{equation}\label{p4}
d_4(\varphi_1,\varphi_2)(\varphi_2')^4-
2d_2(\varphi_1,\varphi_2)(\varphi_2')^2+d_0(\varphi_1,\varphi_2)=0,
\end{equation}
where $d_0$, $d_2$ and $d_4$ are polynomials in $\varphi_1$,
$\varphi_2$ and $p$. The discriminant of (\ref{p4}) is given by
$$\Delta:=-64\varphi_1^2\varphi_2^2w_1w_2w_3,$$
with
$$\begin{array}{l}w_1(\varphi_1,\varphi_2)=\varphi_1^2+\varphi_2^2-1,\\
\\
w_2(\varphi_1,\varphi_2)=p^2\varphi_1^2-(3-4p^2)\varphi_2^2+p^2(3-4p^2),\\
\\
w_3(\varphi_1,\varphi_2)=-(3-4p^2)\varphi_1^2+16p^2\varphi_2^2-4p^2(3-4p^2).
\end{array}$$
It is quite easy to show that, for any $p$, each one of the curves
$(w_i=0)$ contains the orbit of a particular solution of
(\ref{first}). Moreover, the unit circle $(w_1=0)$ represents the
orbit of the solution of (\ref{first}) satisfying the initial
conditions (\ref{datafirst}) with $p=1$. For $p=\sqrt{3/8}$, we
have $w_3\equiv -4w_2$ and the curve $(w_2=0)$ contains the orbit
of the solution of (\ref{first})-(\ref{datafirst}).

These particular algebraic orbits suggest us searching solutions
$(\varphi_1,\varphi_2)$ defined by algebraic relations of the form
$w_4(\varphi_1,\varphi_2)=F(\varphi_1^2,\varphi_2^2)=0$, where $F$
is a polynomial of degree $\le4$. Apart the three quadrics above,
the only additional solution of this type we found is
$$w_4=(\varphi_1^2+4\varphi_2^2)^2-12\varphi_2^2=0.$$
Like $(w_1=0)$, the curve $(w_4=0)$ is independent of $p$ and
represents the orbit of a particular solution of (\ref{first}) for
arbitrary values of $p$. Since
$$w_4(\varphi_1,\varphi_2)=
(\varphi_1^2+4\varphi_2^2-2\sqrt3\,\varphi_2)(\varphi_1^2+4\varphi_2^2+2\sqrt3\,\varphi_2),$$
the set $(w_4=0)$ is the union of two ellipses passing through the
origin, each one being symmetric to the other with respect to the
$\varphi_1$-axis. The upper ellipse
\begin{equation}\label{ellipse}\varphi_1^2+4\varphi_2^2-2\sqrt3\,\varphi_2=0\end{equation}
corresponds to the orbit of the solution of
(\ref{first})-(\ref{datafirst}) associated with
$p=\frac{\sqrt3}2$.

\subsection{Proof of Proposition \ref{pro1}(2): case $\mathbf{ p=\frac{\sqrt3}2}$.}
In this case, the orbit of the solution of
(\ref{first})-(\ref{datafirst}) is given by (\ref{ellipse}). The
only critical point of (\ref{first}) lying on this ellipse is the
origin, which is also a critical point of the system (\ref{s1}).
Therefore, $(\varphi_1(y),\varphi_2(y))$ tends to the origin as
$y$ goes to infinity, and, hence, it is not periodic (see also
\cite{JNP}).

\emph{From now on, we will assume that $p\not=\frac{\sqrt3}2$}.

\subsection{A bounded region for the orbit}

The orbit of the solution of (\ref{first})-(\ref{datafirst}) must
lie in the region of the $(\varphi_1,\varphi_2)$-plane where the
discriminant $\Delta$ of  (\ref{p4}) is nonnegative. This region,
$(\Delta\ge0)$, is a bounded domain delimited by the unit circle
$(w_1=0)$ and the quadrics $(w_2=0)$ and $(w_3=0)$. Its shape
depends on the values of $p$.
\begin{itemize}
\item For $p\in(0,\sqrt3/2)$, $(w_2=0)$ and $(w_3=0)$ are
hyperbolae.

\item The case $p=\sqrt{3/8}$ is a special one since then,
$w_3\equiv-4w_2$, and the region $(\Delta\ge0)$ shrinks to the arc
of the hyperbola $(w_2=0)$ lying inside the unit disk.

\item For $p\in(\sqrt3/2,1]$, $w_3$ is positive and $(w_2=0)$ is
an ellipse.

\end{itemize}

From (\ref{p4}) and (\ref{first}) one can express $\varphi_1'$ and
$\varphi_2'$ in terms of $\varphi_1$, $\varphi_2$ and $p$. Thus,
we obtain a multi-valued 2-dimensional dynamical system
parameterized by $p$ with the initial conditions $\varphi_1(0)=0$
and $\varphi_2(0)=p$. However, the dynamics of such a multi-valued
system is very complex to study. Fortunately, as we will see in
the next subsections, the system (\ref{s1}) can be transformed, by
means of a suitable change of variables, into a Hamiltonian
system, completely integrable by quadratures.

\subsection{Hamiltonian dynamical system }

Let us introduce the new variables $q_1$ and $q_2$ defined by
\begin{equation}\label{q12}
q_1:=\frac1{\sqrt2}\,\varphi_1\hskip 1cm,\hskip1cm
q_2:=\sqrt2\,\varphi_2.\end{equation} The system (\ref{s1})
becomes
\begin{equation}\label{s1q}\left\{\begin{array}{l}
q_1''=[1-4(q_1^2+q_2^2)]q_1=-\frac{\partial V}{\partial q_1},\\ \\
q_2''=4[1-q_1^2-q_2^2]q_2=-\frac{\partial V}{\partial q_2},
\end{array}\right.\end{equation}
with
$$V(q_1,q_2):=(q_1^2+q_2^2)^2-\frac12q_1^2-2q_2^2.$$
Therefore, one has a Hamiltonian system with two degrees of
freedom. The Hamiltonian $H$ is given by
$$H(q_1,q_2,q_1',q_2'):=\frac12[(q_1')^2+(q_2')^2]+V(q_1,q_2).$$
This Hamiltonian is a first integral of (\ref{s1q}) (notice that
$H=\frac14 H_1$). A second independent first integral of
(\ref{s1q}) can be obtained from $H_2$. Consequently, the
Hamiltonian system (\ref{s1q}) is integrable and all its bounded
orbits in phase space $(q_1,q_2,q_1', q_2')$ are contained in a
2-dimensional topological torus (see \cite{Arnold}), which means
that the corresponding solutions are periodic or quasi-periodic,
provided that there is no critical point in the closure of the
orbit. However, it is in general difficult to decide whether such
a solution is periodic or not. The corresponding topological torus
obtained from (\ref{first}) is given by:
\begin{equation}\label{firstq12}\left\{\begin{array}{l}
\frac12[(q_1')^2+(q_2')^2]+(q_1^2+q_2^2)^2-\frac12q_1^2-2q_2^2+p^2(3-4p^2)=0,\\
\\
3q_2^2(q_2^2-2)+3q_1^2q_2^2+(q_1')^2q_2^2-2q_1q_1'q_2q_2'\\ \\
\hskip 3cm +\frac12(3+2q_1^2)(q_2')^2+4p^2(3-4p^2)=0.
\end{array}\right.\end{equation}

It is important to notice that the second first integral is also
quadratic in the $q_1'$ and $q_2'$ variables. Indeed, this enables
us to apply the Bertrand-Darboux-Whittaker theorem: \textit{Given
a Hamiltonian system defined by
$$H=\frac12[(q_1)'^2+(q_2')^2]+V(q_1,q_2),$$
the system admits an additional independent first integral,
quadratic in $q_1'$ and $q_2'$, if and only if the system is
separable in cartesian, polar, parabolic, or elliptic-hyperbolic
coordinates.} (see \cite{AP,W} for details).

In our case, an adequate change of variables is a parabolic one,
given by
\begin{equation}\label{parab}
\left\{\begin{array}{l} q_1^2=-\frac23uv,\\ \\
q_2^2=\frac16(3+2u)(3+2v).
\end{array}\right.\end{equation}
Indeed, from (\ref{firstq12}), one obtains after an elementary
computation
\begin{equation}\label{suv}
\left\{\begin{array}{l} (u')^2 =\frac{P(u)}{(u-v)^2},\\ \\
(v')^2 =\frac{P(v)}{(u-v)^2},\end{array} \right.
\end{equation}
where $P(s):=s(1-2s)(3+2s)(2p^2+s)(3-4p^2+2s)$. Observe that
(\ref{suv}) is not completely decoupled yet; this can be done by
means of a suitable change of the independent variable (see
Subsection 3.6). Each one of the quadrics $(w_1=0)$, $(w_2=0)$ and
$(w_3=0)$ is transformed into two parallel lines. Indeed, we have
$w_1(u,v)=-\frac 14(1-2u)(1-2v)$, $w_2(u,v)=-(\frac
32-2p^2+u)(\frac 32-2p^2+v)$ and $w_3(u,v)=4(2p^2+u)(2p^2+v)$.
Also, we have $\Delta=\frac{16}9P(u)P(v)$. Thus, the region
$(\Delta\ge 0)$ is
transformed into the region $(P(u)P(v)\ge 0)$.\\
Observe that the system (\ref{suv}) is symmetric in $u$ and $v$.
As the change of variables (\ref{parab}) is also symmetric in $u$
and $v$, and since $uv$ must be non positive, one can assume,
without loss of generality, that $u\ge0$ and, hence, $-\frac32\le
v\le0$. Now, the condition $(\Delta\ge0)$ implies that $(u,v)\in
I_1\times I_2$, where
\begin{equation}\label{i1}
I_1:=[\alpha_0,1/2]:=\left\{\begin{array}{ll}
\left[0,\frac 12\right]&\mbox{if }p^2<\frac34,\\
\\
\left[2p^2-\frac 32,\frac12\right]&\mbox{if }\frac34<p^2\le1,\\
\end{array}\right.
\end{equation}
and
\begin{equation}\label{i2}
I_2:=[a_0,a_1]:=\left\{\begin{array}{ll}
\left[2p^2-\frac 32,-2p^2\right]&\mbox{if }p^2\le\frac38,\\
\\
\left[-2p^2,2p^2-\frac 32\right]&\mbox{if }\frac38\le p^2<\frac34,\\
\\
\left[-\frac32,0\right]&\mbox{for }\frac34<p^2\le1.\\
\end{array}\right.
\end{equation}
The initial conditions (\ref{datafirst}) become
$$(u(0), v(0))=\left\{\begin{array}{ll}
(0, 2p^2-\frac 32) &\mbox{if }p^2<\frac34,\\
\\
(2p^2-\frac 32, 0) &\mbox{if }\frac34<p^2\le1.\\
\end{array}\right. $$
In all cases, $u(0)$ and $v(0)$ are zeros of $P$. Hence
$$(u'(0), v'(0))=(0,0).$$
The behavior of $(u,v)$ near $y=0$ is then determined by the
acceleration vector
$$(u''(0), v''(0))=\left\{\begin{array}{ll}
(\frac{12p^2}{3-4p^2}, \frac{16p^2(1 - p^2)(3 - 8p^2)}{(3 - 4p^2)}) &\mbox{if }p^2<\frac34,\\
\\
(\frac{16p^2(1
- p^2)(3 - 8p^2)}{(3 - 4p^2)}, \frac{12p^2}{3-4p^2}) &\mbox{if }\frac34<p^2\le1.\\
\end{array}\right. $$

Notice that for $p=\sqrt\frac38$, $v$ is constant, namely
$v(y)=-\frac34$ for all $y$, while for $p=1$, $u$ is constant with
$u(y) = \frac12$ for all $y$.

\subsection{Proof of Proposition \ref{pro1}(1):
Decoupling the system }

In order to completely decouple the previous system we introduce a
change of the independent variable $y\mapsto \tau$ defined by:
$$\frac{d\tau}{dy}=\frac 1{u-v}.$$
Notice that this change of variable is one-to-one since $u-v\ne0$
(indeed, $I_1\cap I_2=\emptyset$). In this new variable, the
system splits into two independent equations:
\begin{equation}\label{su}
(\dot{u})^2 =P(u),
\end{equation}
\begin{equation}\label{sv}
(\dot{v})^2 =P(v),
\end{equation}
where $\dot{u}:=du/d\tau$ and $\dot{v}:=dv/d\tau$. The solution
$\tau\mapsto u(\tau)$ of (\ref{su}) is also a solution of the
second order ODE
\begin{equation}\label{zpp}\ddot{u}=\frac12P'(u),\end{equation} with the
initial conditions $u(0)=\alpha_0$ (see (\ref{i1}) for the
definition of $\alpha_0$) and $\dot{u}(0)=0$, where $P':=dP/du$.
This solution lies on the curve
\begin{equation}\label{ODE}(\dot{u})^2-P(u)=0\end{equation}
in the $(u,\dot{u})$-phase plane of (\ref{zpp}). Since $\alpha_0$
and $\frac12$ are two consecutive zeros of $P$, the equation
(\ref{ODE}) in the region $\alpha_0\le u\le \frac12$ represents a
closed curve. On the other hand, it is easy to check that $P$ and
$P'$ admit no common zero in the interval $[\alpha_0,\frac12]$.
Hence, there exists no critical point for (\ref{zpp}) on the orbit
defined by (\ref{ODE}) inside the region $\alpha_0\le u\le
\frac12$. Consequently, this closed orbit corresponds to a
periodic solution of (\ref{zpp}) and, therefore, $\tau\mapsto
u(\tau)$ oscillates between $\alpha_0$ and $\frac12$.

A similar analysis for $\tau\mapsto v(\tau)$ implies that it is a
periodic solution of
\begin{equation}\label{vpp}\ddot{v}=\frac12P'(v),\end{equation} with the
initial conditions $v(0)=a_0$ (see (\ref{i2}) for the definition
of $a_0$) and $\dot{v}(0)=0$. Consequently, $\tau\mapsto v(\tau)$
oscillates between $a_0$ and $a_1$. This proves Assertion (1) of
Proposition \ref{pro1}(1).

\subsection{Proof of Proposition \ref{pro1}(3): case $\mathbf{p>{\sqrt3\over 2}}$}

We have just seen that $v(\mathbb{R}) = [a_0,a_1]$, with
$a_0=-\frac32$ for $p\in(\frac{\sqrt3}2,1]$ (see (\ref{i2})). This
implies that $q_2$, and then $\varphi_2$, vanishes at least once
in a period (see (\ref{parab})).

\subsection{About the periods of $u$ and $v$: case $\mathbf{p<{\sqrt3\over 2}}$}

Let us denote $\mathcal{T}_u(p)$ the period of $u$. The function
$\tau\mapsto u(\tau)$ oscillates between $\alpha_0=0$ and
$\frac12$ with velocity $(\dot{u})^2 =P(u)\neq 0$ if $u\in
(0,\frac12)$.  Hence, $u(\tau)$ increases from $0$ to $\frac12$
when $\tau$ goes from 0 to $\mathcal{T}_u(p)\over 2$. It follows
that
\begin{equation}\label{Tu}
\mathcal{T}_u(p)=2\int_0^{\frac12}\frac{ds}{\sqrt{P(s)}}.
\end{equation}
Similarly, for $p\ne \sqrt{3/8}$, the period $\mathcal{T}_v(p)$ of
$\tau\mapsto v(\tau)$ is given by
$$
\mathcal{T}_v(p)=2\int_{a_0}^{a_1}\frac{ds}{\sqrt{P(s)}}.
$$
 Setting, for $p\ne \sqrt{3/8}$, $s=(3-8p^2)r-\frac 32+2p^2$,
one can write
\begin{equation}\label{Tv}
\mathcal{T}_v(p)=2\int_{0}^{\frac 12}\frac{dr}{\sqrt{Q(r)}},
\end{equation}
where
$$Q(r):=2r(1-2r)[2p^2+(3-8p^2)r][2-2p^2-(3-8p^2)r][3-4p^2-2(3-8p^2)r].$$
Hence, the functions $\mathcal{T}_u(p)$ and $\mathcal{T}_v(p)$ are
explicitly given by complete hyper-elliptic integrals. Although
the function $\mathcal{T}_v$ is not defined at $p=\sqrt{3/8}$, its
limit exists. Indeed, setting $t=\frac34+r$ and $\alpha:=\frac
34-2p^2$, we get
$$\mathcal{T}_v(p)=2\int_{-\alpha}^{\alpha}\frac{dt}
{\sqrt{(\frac
32-2t)(\frac52-2t)(\frac32+2t)}\sqrt{\alpha^2-t^2}}.$$ As
$\alpha\to0$, we have
$$
\mathcal{T}_v(p)\sim2\int_{-\alpha}^\alpha\frac{4dt}
{3\sqrt{10}\sqrt{\alpha^2-t^2}}.
$$
Thus
$$\lim_{p\to\sqrt{\frac38}}\mathcal{T}_v(p)=\frac{8\pi}{3\sqrt{10}}.$$
On the other hand, we have
$$\mathcal{T}_u(\sqrt{3/8})=\frac45 \Pi(2/5,1/4),$$
where $\Pi$ is the complete elliptic integral of the third kind
given by
$$\Pi(n,m):=\int_0^{\frac\pi2}\frac{d\theta}{(1-n\sin^2\theta)\sqrt{1-m\sin^2\theta}}.$$
Since for $p=\sqrt{3/8}$, $v$ is constant, the
couple $(u,v)$ is periodic of period $\mathcal{T}_u(\sqrt{3/8})$.\\

\noindent\underline{Behavior of ${\mathcal{T}_v-\mathcal{T}_u}$
and $ \mathcal{T}_v/\mathcal{T}_u$ near ${p=0}$}: One has
$$\mathcal{T}_v(p)-\mathcal{T}_u(p)=2\int_0^{\frac 12}\left[
\frac1{\sqrt{Q(s)}}- \frac1{\sqrt{P(s)}}\right]\,ds.$$ The
integral of $ \frac1{\sqrt{P(s)}}$ is singular only at $p=0$. A
direct computation gives
$$\int_0^{\frac 12}\frac{ds}{\sqrt{P(s)}}\sim \int_0^{\frac
12}\frac{ds}{3\sqrt{s^2+2p^2s}}\sim-\frac 23\ln(p).$$ Similarly,
we get
$$\int_0^{\frac 12}\frac{ds}{\sqrt{Q(s)}}\sim -\ln p.$$
In other words $\mathcal{T}_v(p)-\mathcal{T}_u(p)\to+\infty$ as
$p\to0$ while the ratio $ \mathcal{T}_v(p)/\mathcal{T}_u(p)$ goes
to $\frac 32$ (see the figure below).\\

\noindent\underline{Behavior of ${\mathcal{T}_v-\mathcal{T}_u}$
and $ \mathcal{T}_v/\mathcal{T}_u$ near ${p=\sqrt3/2}$}: As
$p\to\sqrt3/2$, $\mathcal{T}_u(p)\sim-\frac23\ln(\sqrt3/2-p)$ and
$\mathcal{T}_v(p)\sim-\ln(\sqrt3/2-p)$. Hence,
$\mathcal{T}_v-\mathcal{T}_u\to+\infty$ as $p\to\sqrt3/2$ and
$ \mathcal{T}_v(p)/\mathcal{T}_u(p)$ goes again to $\frac32$.\\

Thus, $\mathcal{T}_v(p)>\mathcal{T}_u(p)$ near $p=0$ and
$p=\sqrt3/2$. Actually, one has
$\mathcal{T}_v(p)>\mathcal{T}_u(p)$ for all $p\in (0,\sqrt3/2)$ as
shown by the graphic representation of $\mathcal{T}_v$ and
$\mathcal{T}_u$ given below.

\begin{center}
\includegraphics[angle=0,width=8cm]{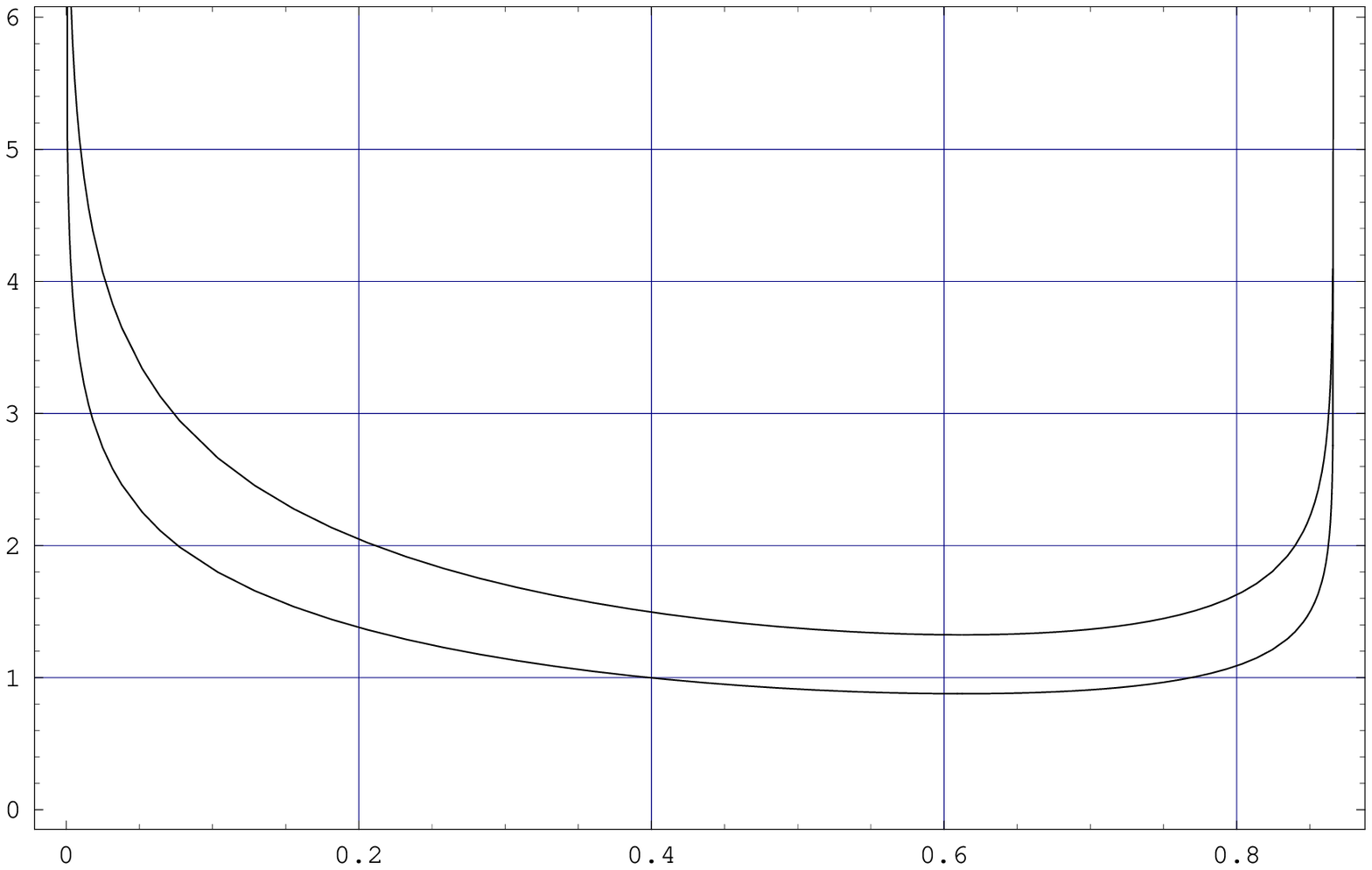}

{\small The functions $p\mapsto \mathcal{T}_v(p)$ (the upper one)
and $p\mapsto\mathcal{T}_u(p)$.}
\end{center}

\subsection{Proof of Proposition \ref{pro1}(4).}
The couple $(u,v)$ is periodic if and only if the ratio
$R(p):=\frac{\mathcal{T}_v(p)}{\mathcal{T}_u(p)}$ is a rational
number. From the previous subsection, $R$ is a nonconstant
continuous function on $(0,\sqrt3/2)$ with $\lim_{p\to
0}R(p)=\lim_{p\to \sqrt3/2}R(p)= 3/2$. The range of $R$ is a
closed interval $[r_1,r_2]\subset[1.480473,1.507784]$ (see the
figure below).

\begin{center}
\includegraphics[angle=0,width=8cm]{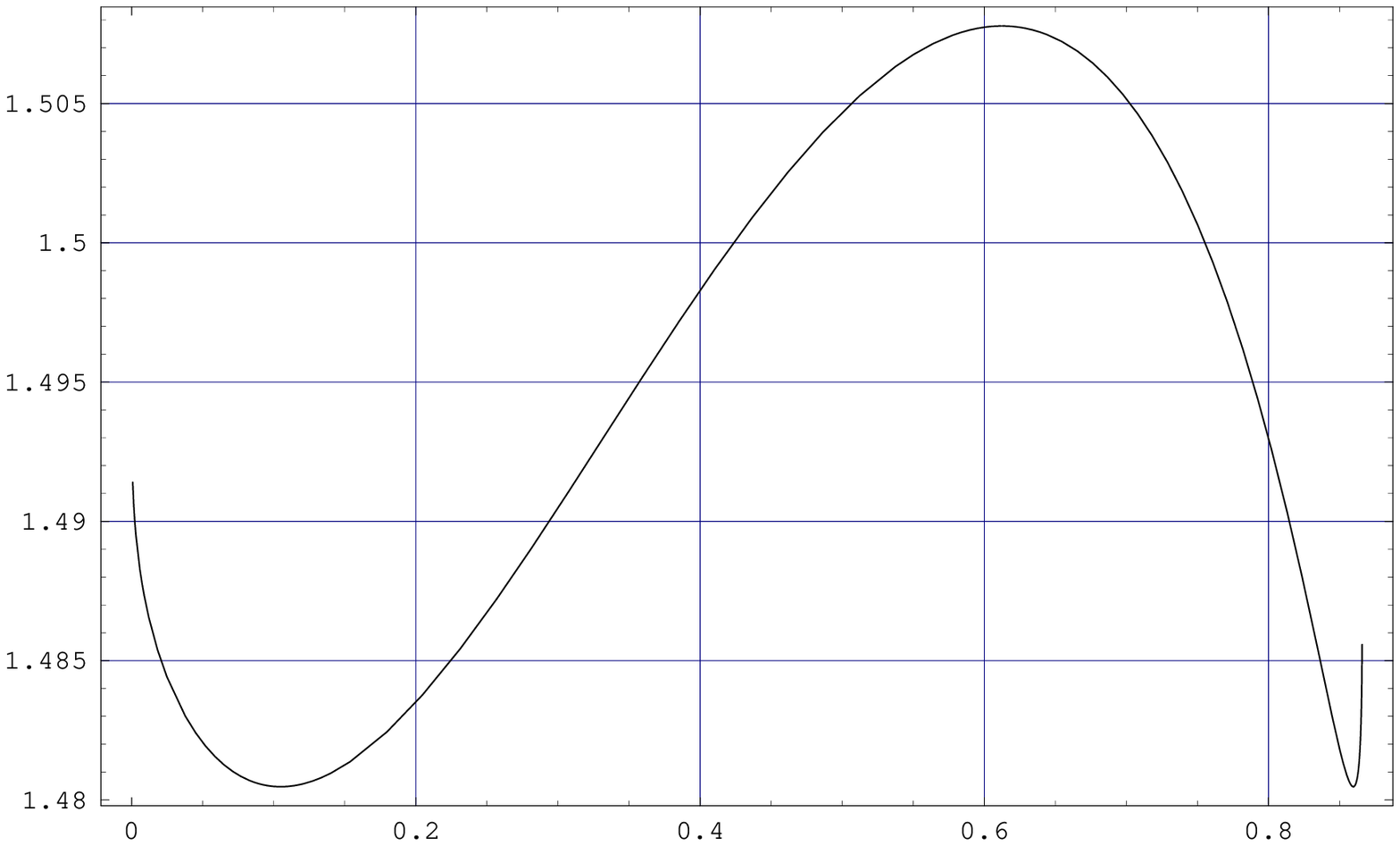}

{\small The ratio $\frac{\mathcal{T}_v}{\mathcal{T}_u}$}
\end{center}

To end the proof of Assertion (4), we only need to define
$\mathcal{P}$ to be the set of $p\in(0,\sqrt3/2)$ such that $R(p)$
is a rational number.

\subsection{Proof of Proposition \ref{pro1}(5).}
In the case $p=\sqrt{3/8}$ we have, $\forall y\in\mathbb{R}$,
$v(y)=-\frac34$ and, then, $\varphi_1^2=u$ and $\varphi_2^2
=\frac18 (3+2u)$. The couple of periodic functions $(\varphi_1,
\varphi_2)=(\sqrt u ,\sqrt{\frac18 (3+2u)})$ on
$[0,\mathcal{T}_u(\sqrt{3/8})]$, such that $\varphi_1$ is odd and
$\varphi_2$ is even, solves the original system and satisfies
Condition (\ref{zeros}).

Let $p\in \mathcal{P}$, $p\neq\sqrt{3/8}$, and let $\frac qm \in
\mathbb{Q}$ be an irreducible fraction, with $q$, $m \in
\mathbb{N}$, such that
$R(p)=\frac{\mathcal{T}_v(p)}{\mathcal{T}_u(p)}=\frac qm$. The
period of the couple $(u,v)$ is given by
$$\mathcal{T}(p)=q\mathcal{T}_u(p)=m\mathcal{T}_v(p).$$
The number of zeros of $u$ in a period, for instance $
[0,\mathcal{T}(p))$, of $(u,v)$ is equal to $q$ times the number
of zeros of $u$ in $[0,\mathcal{T}_u(p))$. As we saw above,
$\forall p\in(0,\sqrt3/2)$, one has $1<R(p)<2$. Hence, $m\ge 2$
and $q>m$, which implies $q\ge3$. Since $u(0)=0$, the number of
zeros of $u$ in a period of $(u,v)$ is at least 3. Since
$\varphi_1$ is odd, the period of $(\varphi_1,\varphi_2)$ is twice
the period of $(u,v)$. From  $\varphi_1= 2\sqrt{-\frac13 uv}$ on
$[0,\mathcal{T}(p))$, one deduces that
 $\varphi_1$ admits at least 6 zeros in a period of
$(\varphi_1,\varphi_2)$. Notice that the case $p=\sqrt{3/8}$ is
special since, for this value of $p$, $v$ is constant and the
couple $(u,v)$ is periodic whose period is equal to that of $u$.

\subsection{On the shape of solutions}

Although it is not necessary for the proof of our main result, one
can obtain as a by product of our study, some properties
concerning the shape of solutions. First, notice that, for
$p\in(0,\frac{\sqrt3}2)$, the 2-dimensional dynamical system
admits four critical points in the region
$(\Delta\ge0)\cap(\varphi_2\ge0)$: $A=(2p/\sqrt3,
\sqrt{1-4p^2/3})$, $B=( \sqrt{1-4p^2/3} , 2p/\sqrt3)$ and their
symmetric with respect to the $\varphi_2$-axis that we denote $A'$
and $B'$. Notice that these critical points are on the boundary of
the region $(\Delta\ge0)$.

\noindent\underline{Non-periodic solutions.} They correspond to
the case where $\mathcal{T}_v(p)/\mathcal{T}_u(p)$ is irrational.
In this case, the orbit $(u,v)$ fill the rectangle $I_1\times
I_2$, and, then, in the $(\varphi_1,\varphi_2)$-plane, the orbit
fill the region $(\Delta\ge0)$.

\noindent\underline{Periodic solutions.}  For $p=\sqrt{3/8}$, the
solution $(\varphi_1,\varphi_2)$ lies on the hyperbola of equation
$\varphi_1^2-4\varphi_2^2+3/2=0,$ oscillating between the points
$A=(-\frac1{\sqrt2},\frac1{\sqrt2})$ and
$A'=(\frac1{\sqrt2},\frac1{\sqrt2})$. For $p\neq\sqrt{3/8}$, let
$q/m=\mathcal{T}_v(p)/\mathcal{T}_u(p)$ be an irreducible
fraction, with $q,m\in\mathbb{N}$, and set
$\mathcal{T}:=q\mathcal{T}_u(p)=m\mathcal{T}_v(p)$. Geometrically,
this means that $u$ makes $q$ round trips in a period while $v$
makes $m$ round trips. We distinguish three cases:
\begin{itemize} \item If $q$ and $m$ are both odd, then
$u(\mathcal{T}/2)=\frac12$ and $v(\mathcal{T}/2)=a_1$. This
corresponds to the point $A$ for $p^2<\frac38$ and to the point
$B$ for $\frac 38<p^2<\frac34$. The orbit is not closed and
$(\varphi_1,\varphi_2)$ oscillates between $A$ and $A'$ or $B$ and
$B'$.

\item If $q$ is odd and $m$ is even, then
$u(\mathcal{T}/2)=\frac12$ and $v(\mathcal{T}/2)=a_0$. This
corresponds to the point $B$ for $p^2<\frac38$ and to the point
$A$ for $\frac 38<p^2<\frac34$. Again, $(\varphi_1,\varphi_2)$
oscillates between $A$ and $A'$ or $B$ and $B'$.

\item If $q$ is even and $m$ is odd, then  $u(\mathcal{T}/2)=0$
and $v(\mathcal{T}/2)=a_1$. This corresponds in the
$(\varphi_1,\varphi_2)$-plane to the point $(0,\frac 34-p^2)$.
This point is the intersection between the quadric $(w_3=0)$ with
the $\varphi_2$-axis. In this case the orbit is closed.
\end{itemize}

\end{document}